\documentstyle[12pt]{article}
\begin{document}
\def\dsp{\displaystyle}
\def\Rr{{bf R}}
\def\Zz{bf Z}
\def\Nn{bf N}
\def\get{\hbox{{\goth g}$^*$}}
\def\g{\gamma}
\def\om{\omega}
\def\r{\rho}
\def\a{\alpha}
\def\s{\sigma}
\def\vfi{\varphi}
\def\l{\lambda}
\def\implique{\Rightarrow}
\def\o{{\circ}}
\def\Diff{\hbox{\rm Diff}}
\def\S1{\hbox{\rm S$^1$}}
\def\Hom{\hbox{\rm Hom}}
\def\Vect{\hbox{\rm Vect}}
\def\const{\hbox{\rm const}}
\def\ad{\hbox{\hbox{\rm ad}}}
\def\semid{\hbox{\bb o}}
\def\blanc{\hbox{\ \ }}

\def\pds#1,#2{\langle #1\mid #2\rangle} 
\def\f#1,#2,#3{#1\colon#2\to#3} 

\def\hfl#1{{\buildrel{#1}\over{\hbox to
12mm{\rightarrowfill}}}}

\title{Sturm theory, Ghys theorem on zeroes of the Schwarzian derivative
and flattening of Legendrian curves}

\author{
V. Ovsienko\\
{\small CNRS, Centre de Physique Th\'eorique}
\thanks{CPT-CNRS, Luminy Case 907,
F--13288 Marseille, Cedex 9
FRANCE}
\and
S. Tabachnikov\thanks{Supported in part by an NSF grant DMS-9402732}\\
{\small University of Arkansas
\thanks{UARK, Department of Mathematics, Fayetteville, AR 72701 USA} and}\\
{\small Wolfson College, University of Cambridge}
}

\date{}

\maketitle

\thispagestyle{empty}

Etienne Ghys has recently discovered a beautiful theorem:
{\it given a diffeomorphism of the projective line, there exist
at least four distinct points in which the diffeomorphism is
unusually well approximated by projective transformations} - \cite{ghy}.
The points in question are the ones in which the 3-jet of the diffeomorphism
is that of a projective transformation;
in a generic point the order of approximation is only 2.
In other words, the Schwarzian derivative
of every diffeomorphism of ${\bf RP}^1$ has at least four distinct zeroes.

The theorem of Ghys is analogous to the classical four vertex theorem:
a closed convex plane curve has at least four curvature extrema \cite{muh}.
The proof presented by E. Ghys was purely geometrical: it was inspired by the
Kneser theorem on osculating circles of a plane curve \cite{kne}.

We will prove an amazing strengthening
of the classical Sturm comparison theorem
and deduce Ghys' theorem from it.

Sturm theory is related to the geometry of curves (see
\cite{arn1},\cite{arn2},\cite{arn3},\cite{arn4},\cite{gmo},\cite{tab1},\break\cite{tab2}).
This relation is based on the following general idea.
A geometrical problem is reformulated as the problem on the least number of
zeroes of
a certain function; this function turnes to be orthogonal to
a number of functions enjoing a disconjugacy property.
This implies an estimate below on the number of zeroes.

An example of such an estimate is the Hurwitz theorem:
{\it the number of zeroes of a function on the circle
is minorated by that of its first harmonic} \cite{hur}.

We will show that the Ghys theorem is equivalent to the following result.

Consider linear symplectic $({\bf R}^4,\omega)$, choose
two orthogonal (with respect to the symplectic structure)
symplectic subspaces ${\bf R}^2_1$ and ${\bf R}^2_2$.
 Denote by ${\bf RP}^1_1$ and ${\bf RP}^1_2$
their projectivizations in the contact ${\bf RP}^3$.

\proclaim Theorem 1. Let
$C\subset{\bf RP}^3$ be a closed Legendrian curve
such that its projection on ${\bf RP}^1_1$
from ${\bf RP}^1_2$ and its projection on ${\bf RP}^1_2$
from ${\bf RP}^1_1$
are diffeomorphisms. Then
$C$ has at least 4 flattening points.\par

We will show that this is a reformulation of the Ghys theorem.

\section{Symplectization of the projective line \break
diffeomorphism
and the Sturm-Liouville  equation}
The projective line ${\bf RP}^1$ is the space of lines
through the origin in the plane ${\bf R}^2$.
A linear transformation of the plane induces a transformation of ${\bf RP}^1$
called a projective transformation.
Projective transformations form a group $PSL_2$.

Given a diffeomorphism
$
f:{\bf RP}^1\rightarrow {\bf RP}^1
$
and a point $x$ of ${\bf RP}^1$, there exists a unique projective
transformation
$g\in PSL_2$ that approximates $f$ in $x$ up to the second order. That is,
$j^2_x(f)=j^2_x(g)$. We are interested in the points in which
$j^3_x(f)=j^3_x(g)$. Call such points {\it projective points}
 of the diffeomorphism $f$.
Without loss of generality assume that $f$ preserves the orientation.

Given an orientation preserving diffeomorphism $f:{\bf RP}^1\rightarrow {\bf
RP}^1$,
there exists a unique area preserving homogeneous (of degree 1) diffeomorphism
$F$ of
the punctured plane ${\bf R}^2\setminus \{0\}$
that projects to $f$.
Let $\alpha$ be the angular parameter on
${\bf RP}^1$ so that $\alpha$ and $\alpha+\pi$ correspond to the same point.
Denote by $f(\alpha)$ the lift of the diffeormorphism $f$
to a diffeomorphism of the line such that
$f(\alpha+\pi)=f(\alpha)+\pi$.
 The diffeomorphism $F$ is given in polar coordinates by
$$
F:(\alpha,r)\mapsto(f(\alpha),\;r\dot f ^{-1/2}(\alpha))
$$
where the dot denotes $d/d\alpha$.
Projective points of the diffeomorphism $f$ give rise to lines along which
the symplectomorphism $F$ is unusualy well (up to the second jet)
approximated by a linear
area preserving transformation of ${\bf R}^2$.

Let $\gamma(\alpha)$ be the image of the unit circle under $F$.
 It is
a centrally symmetric curve: $\gamma(\alpha+\pi)=-\gamma(\alpha)$, and it
 bounds area $\pi$.
If $F\in SL_2$ then the corresponding curve is a central ellipse
of area $\pi$.
Thus for every point $x\in \gamma$ there exists a unique central ellipse
of area $\pi$ tangent to $\gamma$ in $x$.
The projective points
of the diffeomorphism $f$ correspond to the points of $\gamma$
in which the tangent ellipse has the second order contact with $\gamma$.

\null

\vspace{6cm}

{}~\hspace{2cm}\special{illustration figghys_1.ps}

\vspace{0,2cm}

\proclaim Lemma 1.1. The parametrised curve $\gamma(\alpha)$ satisfies
the equation:
\begin{equation}
\ddot\gamma(\alpha)=-k(\alpha)\gamma
\label{sl}
\end{equation}
where $k(\alpha)$ is a $\pi$-periodic smooth function.\par

{\bf Proof}. Let $\gamma_0(\alpha)$ be the parametrised unit circle.
Then $[\gamma_0,\dot\gamma_0]=1$ where $[\;,\;]$ is the oriented area
of the parallelogram generated by two vectors.
Applying the symplectomorphism $F$ one obtains $[\gamma,\dot\gamma]=1$.
Differentiate to obtain $[\gamma,\ddot\gamma]=0$. Hence, the vector
$\ddot\gamma$ in collinear to $\gamma$. The result follows.

\proclaim Lemma 1.2. The projective points
of the diffeomorphism $f$ correspond to the points in which
$k(\alpha)=1$.\par

{\bf Proof}. If $F\in SL_2$ the curve $\gamma$ is a central ellipse
of area $\pi$. In this case $k(\alpha)\equiv 1$. Projective points
are points of second order contact with central ellipses of area $\pi$.
Thus $k(\alpha)=1$ in these points.

\vskip 0,3cm

In view of the above lemmas it suffices to prove that $k(\alpha)-1$
has at least 4 distinct zeroes on the interval $[0,\pi]$.

\section{Strengthened Sturm comparison theorem}
The Sturm comparison theorem is stated as follows.
{\it Given two Sturm-Liouville equations
$\ddot\phi(\alpha)=-k(\alpha)\phi(\alpha)$ with potentials
$k_1(\alpha)>k_2(\alpha)$, let $\phi_1$ and $\phi_2$ be
solutions of the respective equations. Then between any two zeros of
$\phi_2$ there exists a zero of $\phi_1$} \cite{stu}.
Equivalently, if $\phi_1$ and $\phi_2$ have two coinciding consecutive zeroes
(see fig.2), then
there is a zero of the function $k_1-k_2$ on this interval.

\null

\vspace{4cm}

{}~\hspace{-0,7cm}\special{illustration figghys_2.ps}

\vspace{0,2cm}

This already implies the existence of two projective points
of a diffeomorphism of ${\bf RP}^1$. Indeed, every solution of the equation
$\ddot\phi(\alpha)=-k(\alpha)\phi(\alpha)$
 is a linear coordinate (e.g. the projection on the vertical axis)
of the curve $\gamma(\alpha)$ satisfying equation (\ref{sl}).
Therefore, the solutions of $\ddot\phi(\alpha)=-k(\alpha)\phi(\alpha)$
 are antiperiodic on $[0,\pi]$
as well as those of the equation $\ddot\phi(\alpha)=-\phi(\alpha)$ (see fig.3).

\vfill\eject

\null

\vspace{4cm}

{}~\hspace{1,3cm}\special{illustration figghys_3.ps}

\vspace{0,2cm}

We prove the following strengthened version of the Sturm theorem which implies
the theorem of Ghys.

Call a Sturm-Liouville equation
with a $\pi$-periodic potential {\it disconjugate}
if for every solution $\phi(\alpha)$ one has $\phi(\alpha +\pi)=-\phi(\alpha)$,
and every
solution has exactly one zero on $[0,\pi)$. Notice that the Sturm-Liouville
equation
corresponding to a diffeomorphism of ${\bf RP}^1$ is disconjugate
(since the curve $\gamma$ in fig.1 is star-shaped).

\proclaim Theorem 2.1. Given two disconjugate
Sturm-Liouville equations
$\ddot\phi(\alpha)=-k_1(\alpha)\phi(\alpha)$ and
$\ddot\phi(\alpha)=-k_2(\alpha)\phi(\alpha)$, there exist at least 4 distinct
zeroes
of the function $k_1-k_2$ on every interval $[\alpha,\alpha+\pi)$.\par

{\bf Proof}. Let $\phi_1$, $\phi_2$ be solutions of the two Sturm-Liouville
equation, respectively. Then,
$$
\int_{\alpha}^{\alpha+\pi}(k_1-k_2)\phi_1\phi_2\;d\alpha=0
$$
Indeed,
$$
\begin{array}{rcl}
\noalign{\smallskip}
{\dsp 0=\int_{\alpha}^{\alpha+\pi}(\phi_1(\ddot\phi_2+k_2\phi_2)-
\phi_2(\ddot\phi_1+k_1\phi_1))\;d\alpha}&=&
{\dsp \left(\phi_1\dot\phi_2-\phi_2\dot\phi_1)
\right|_{\alpha}^{\alpha+\pi}+}\\
 \noalign{\smallskip}
&&{\dsp \;\;\;\;\;\int_{\alpha}^{\alpha+\pi}(k_2-k_1)\phi_1\phi_2\;d\alpha}
\end{array}
$$
The first term at the right hand side vanishes because
$\phi_1\phi_2$ is $\pi$-periodic.

The theorem is a corollary of the following fact.

\proclaim Lemma 2.2. If a $\pi$-periodic
function $\psi$ is orthogonal to the product
of any two solutions $\phi_1$ and $\phi_2$ of two disconjugate
Sturm-Liouville equation, respectively, then $\psi$ has at least 4 distinct
zeroes
on any interval $[\alpha,\alpha+\pi)$.\par

\vfill\eject

\null

\vspace{4,5cm}

{}~\hspace{-0,7cm}\special{illustration figghys_4.ps}

\vspace{0,2cm}

{\bf Proof}. First prove that $\psi$ has 2 zeroes.
If not, choose $\phi_1$ and $\phi_2$ to have the same zero
on $[\alpha,\alpha+\pi)$ (fig.4).
Since both equations are disconjugate,
neither $\phi_1$ nor $\phi_2$ have other zeroes
on the interval. Then, the product
$\psi\phi_1\phi_2$ is either positive almost everywhere
or negative almost everywhere.

Secondly, assume that $\psi$ changes sign only twice at points
$\alpha_1,\alpha_2$. Choose $\phi_1$ to have zero at $\alpha_1$
and $\phi_2$ -- at $\alpha_2$. Again, $\phi_1$ and $\phi_2$
have no extra zeroes. As before, the integral of
$\psi\phi_1\phi_2$ cannot vanish. The lemma is proved.

{\bf Remark}. The above lemma is a particular case of
a hierarchy of theorems estimating below the number of zeroes
of functions orthogonal to products of solutions of disconjugate
linear differential equations -- see \cite{gmo} where such
theorems are used to prove the existence of 6 affine vertices
of a closed convex plane curve.

\section{Zeroes of the Schwarzian derivative}

The Schwarzian derivative measures the failure of a diffeomorphism
of ${\bf RP}^1$ to be projective.

More specifically, let $f$
be a diffeomorphism, $\alpha\in{\bf RP}^1$ be a point.
Consider four ``infinitly close'' points
$\alpha,\;\alpha+\epsilon,\;\alpha+2\epsilon,\;\alpha+3\epsilon$.
Apply $f$ to obtain the points
$f(\alpha),\;f(\alpha+\epsilon),\;f(\alpha+2\epsilon),\;f(\alpha+3\epsilon)$.
Define the Schwarzian derivative $S(f)(\alpha)$ by:
$$
[f(\alpha),f(\alpha+\epsilon),f(\alpha+2\epsilon),f(\alpha+3\epsilon)]=
[\alpha,\alpha+\epsilon,\alpha+2\epsilon,\alpha+3\epsilon]+
\epsilon^2S(f)(\alpha)+O(\epsilon^3)
$$
where $[\;,\;,\;,\;]$ denotes the cross-ratio of four lines through
the origin.

\vskip 1cm

$$\;\;\;\;\;\;\;\;\;\;\;\;\;\;\;\;\;\;
\;\;\;\;\;\;\;\;\;\;\;\;\;\;\;\;\;\;\;\;\;\;\;\;\;\;\;\;\;\;\;\;\;\;\;\;
[a,b,c,d]=$$
$$\;\;\;\;\;\;\;\;\;\;\;\;\;\;\;\;\;\;
\;\;\;\;\;\;\;\;\;\;\;\;\;\;\;\;\;\;\;\;\;\;\;\;\;\;\;\;\;\;\;\;\;\;\;\;
[A,B,C,D]=$$
$$\;\;\;\;\;\;\;\;\;\;\;\;\;\;\;\;\;\;
\;\;\;\;\;\;\;\;\;\;\;\;\;\;\;\;\;\;
\;\;\;\;\;\;\;\;\;\;\;\;\;\;\;\;\;\;\;\;\;\;\;\;\frac{(A-C)(B-D)}{(B-C)(A-D)}$$

\null

\vspace{0,5cm}

{}~\hspace{0cm}\special{illustration figghys_5.ps}

\vspace{0,2cm}

By the very definition {\it projective points of a diffeomorphism
are zeroes of its Schwarzian derivative}.

The following formula is the result of a direct computation.

\proclaim Lemma 3.1. In the angular parameter $\alpha$
\begin{equation}
S(f)=\frac{\buildrel{\dots}\over f}{\dot f}-
\frac{3}{2}\left(\frac{\ddot f}{\dot f}\right)^2+
2(\dot f^2-1)
\label{sc}
\end{equation}\par

Notice that in the affine parameter $x=\hbox{tg}\alpha$
the Schwarzian derivative is given by the standard formula:
$$
S(f)=\frac{f'''}{f'}-
\frac{3}{2}\left(\frac{f''}{f'}\right)^2
$$
where $'=d/dx$.

{\bf Remark}. The Schwarzian derivative is invariantly defined
as the unique 1-cocycle on the group of diffeomorphisms
of ${\bf RP}^1$ whose kernel is $PSL_2$
with values in the space of quadratic differentials.

The relation between the Schwarzian derivative and the potential
of the Sturm-Liouville equation (\ref{sl}) is as follows.

\proclaim Lemma 3.2.
$\;\;\;\;\;\;\;\;\;\;\;\;\;\;\;\;\;\;\;\;\;\;\;\;k=\frac{1}{2}S(f)+1$.\par

{\bf Proof}. The curve $\gamma(\alpha)$
is the image of the circle under the symplectomorphism $F$.
In Cartesian coordinates,
$$
\gamma(\alpha)=(\dot f^{-1/2}(\alpha)\cos f(\alpha),\;\dot f^{-1/2}(\alpha)\sin
f(\alpha))
$$
Differentiate twice to obtain the result.

\vskip 0,3cm

{\bf Remark}. It is interesting to consider the infinitesimal
analog of the condition $S(f)(\alpha)=0$. Let
$f(\alpha)=\alpha+\epsilon a(\alpha)$ where $a(\alpha)$
is a $\pi$-periodic function.
Then, in the first order in $\epsilon$ one gets from (\ref{sc}) the equation:
$$
\buildrel{\dots}\over a(\alpha)+4\dot a(\alpha)=0
$$
The existence of at least 4 roots of this equation on $[0,\pi)$
follows from the Hurwitz theorem.
Indeed, the space of $\pi$-periodic functions has the basis:
$$
1,\sin 2\alpha,\cos 2\alpha,\sin 4\alpha,\dots
$$
and the function $\buildrel{\dots}\over a+4\dot a$
does not contain the first harmonics.

\section{Projective points of diffeomorphisms as
flattening points of Legendrian curves in projective space}

In this section we prove Theorem 1.

Consider the graph of the homogeneous symplectomorphism $F$.
It is a conical Lagrangian submanifold in the symplectic space
${\bf R}^2\times{\bf R}^2$, the symplectic structure being
$\omega_1\ominus\omega_2$ where $\omega_1$ and $\omega_2$
are the symplectic structures on the factors.
Denote by ${\cal F}$ the projectivization of the graph.
${\cal F}$ is a Legendrian curve in contact ${\bf RP}^3$.

\null

\vspace{6cm}

{}~\hspace{0cm}\special{illustration figghys_6.ps}

\vspace{0,2cm}

Space ${\bf RP}^3$ contains two distinguished projective lines,
the projectivizations of the source and target ${\bf R}^2$.
Both projections of ${\cal F}$ from one of these lines to
another are diffeomorphisms.

Projective points of the diffeomorphism $f$ correspond to
lines along which the graph of the symplectomorphism $F$ has contact of second
order with graphs of linear symplectic mappings. In projectivization
projective points of the diffeomorphism $f$ correspond to
{\it inflection points} of the curve ${\cal F}$.

An inflection point of a curve in ${\bf RP}^3$ is a point in which
the acceleration vector is proportional to the velocity vector.
Said differently, the curve has the second order contact with
its tangent line.
A generic curve in ${\bf RP}^3$ does not have inflection points
(since the inflection point condition has codimension 2).

Conversely, let
$C\subset{\bf RP}^3$ be a closed Legendrian curve
such that its projection on ${\bf RP}^1_1$
from ${\bf RP}^1_2$ and its projection on ${\bf RP}^1_2$
from ${\bf RP}^1_1$
are diffeomorphisms. Then the composition of one
projection with the inverse of another is a diffeomorphism
of ${\bf RP}^1$. The inflection points of $C$
are the projective points of this diffeomorphism.
Theorem 1 is a consequence of
the Ghys theorem and Proposition 4.1 below.

To formulate this proposition recall that {\it flattening points}
of a curve $C\subset{\bf RP}^3$
are the points in which the curve has contact of second order
with its osculating plane (i.e., the vectors
$\dot C,\ddot C,\buildrel{\dots}\over C$
are linearly dependent).
Clearly, inflection points are flattening points.
It turns out that for Legendrian curves the two
notions coincide.

\proclaim Proposition 4.1. The flattening points of a Legendrian curve
in the standard contact ${\bf RP}^3$ are its inflection points.\par

{\bf Proof}. Consider a parametrization $C(t)$. We need to prove
that if the vectors $\dot C,\ddot C,\buildrel{\dots}\over C$
are linearly dependent then so are $\dot C,\ddot C$
(which makes sense in any affine chart and does not depend on its choice).
Let $\theta$ be the standard contact form. The flattening condition reads:
$\theta\wedge d\theta(\dot C,\ddot C,\buildrel{\dots}\over C)=0$.

Let $\widetilde C$ be a lift of $C$ to linear symplectic space ${\bf R}^4$,
$\omega$ -- the linear syplectic structure. Since $C$ is Legendrian,
$\omega(\widetilde C,\dot{\widetilde C})=0$. Differentiate to obtain:
$\omega(\widetilde C,\ddot{\widetilde C})=0$. This means that the
acceleration vector $\ddot C$ belongs to the contact plane (cf. \cite{app}).

Differentiate once again:
$$
\omega(\dot{\widetilde C},\ddot{\widetilde C})+
\omega(\widetilde C,\buildrel{\dots}\over {\widetilde C})=0
$$
The flattening condition reads:
$\omega\wedge\omega(\widetilde C, \dot{\widetilde C},\ddot{\widetilde C},
\buildrel{\dots}\over {\widetilde C})=0$.
In view of the two previous formul{\ae}
$$
0=\omega\wedge\omega(\widetilde C, \dot{\widetilde C},\ddot{\widetilde C},
\buildrel{\dots}\over {\widetilde C})=
\omega(\widetilde C, \buildrel{\dots}\over {\widetilde C})
\omega(\dot{\widetilde C},\ddot{\widetilde C})=
-\omega(\dot{\widetilde C},\ddot{\widetilde C})^2
$$

If this is zero then $d\theta(\dot C,\ddot C)^2=0$. Since
$\dot C,\ddot C$ lie in the contact plane and $d\theta$
in nondegenerate therein, the vectors $\dot C,\ddot C$
are linearly dependent. The proposition is proved.

{\bf Remark}. Parametrize
${\bf RP}^1_1$ by the angular parameter $\alpha$.
Then $\alpha$ also parametrizes $C$.
Let $\widetilde C(\alpha)$ be a lift of $C(\alpha)$ to ${\bf R}^4$ such that
its projection to ${\bf R}^2_1$ has coordinates
$(\cos(\alpha),\sin(\alpha))$. Let $(\phi_1(\alpha),(\phi_2(\alpha))$
be the projection of $\widetilde C(\alpha)$ to ${\bf R}^2_2$.
Then
$$
\left| \begin{array}{cc}
\phi_1&\phi_2 \\
\dot\phi_1&\dot\phi_2\\
\end{array}\right|=\hbox{const}
$$
and $\phi_1$ and $\phi_2$ satisfy
a Sturm-Liouville equation: $\ddot\phi+k\phi=0$
where $k(\alpha)$ is a $\pi$-periodic function.

One gets:
\begin{equation}
\omega\wedge\omega(\widetilde C, \dot{\widetilde C},\ddot{\widetilde C},
\buildrel{\dots}\over {\widetilde C})=\left| \begin{array}{cc}
\phi_1+\ddot\phi_1&\phi_2+\ddot\phi_2 \\
\dot\phi_1+\buildrel{\dots}\over\phi_1&\dot\phi_2+\buildrel{\dots}\over\phi_2\\
\end{array}\right|
\label{fl}
\end{equation}
This in turn is equal to
$\hbox{const}(k(\alpha)-1)^2$ (cf. Theorem 2.1).

\vskip 0,3cm

The Legendrian condition on the curve $C$ in Theorem 1 was essential.

\vskip 0,3cm

{\bf Example}. Consider the curve
$\widetilde C(\alpha)\subset{\bf R}^4$ whose projection
to ${\bf R}^2_1$ is $(\cos \alpha, \sin \alpha)$,
and to ${\bf R}^2_2$ is
$(\cos \alpha+\epsilon\cos 3\alpha,\; \sin \alpha+\epsilon\sin 3\alpha)$.
Then for sufficiently small $\epsilon$ both projections
of $C$ to ${\bf RP}^1_1$ and ${\bf RP}^1_2$ are
diffeomorphisms. However, $C$ has no flattening points:
the determinant (\ref{fl}) equals $192\epsilon^2$.

\section{Inflections of the characteristic curve of a projective line
diffeomorphism}

Consider the Hopf fibration $\pi: RP^3 \rightarrow S^2$ whose fibers are the
characteristic
 lines of the standard contact form in ${\bf {\bf RP}}^3$. Consider also the
fibration
$\widetilde \pi: {\bf R}^2_1 \times {\bf R}^2_2\setminus \{0\} \rightarrow S^2$
that covers $\pi$.

Given a Legendrian curve in ${\bf RP}^3$ its projection is an immersed smooth
curve on
$S^2$ that bounds the area which is a multiple of half that of the sphere. The
Legendrian lines project to great circles in $S^2$ (a one-parameter family of
Legendrian lines over each great circle). The points $\widetilde \pi ({\bf
R}^2_1)=
\pi ({\bf RP}^1_1)$ and $\widetilde \pi ({\bf R}^2_2) = \pi ({\bf RP}^1_2)$
 are antipodal; we think of
them as the poles of the sphere.

We have associated  the Legendrian curve ${\cal F} \subset {\bf RP}^3$ with an
orientation
preserving diffeomorphism of the projective line.

{\bf Definition}. The characteristic curve $c_f$ of a diffeomorphism $f$ is the
curve
$\pi ({\cal F})$.

\proclaim Lemma 5.1. The projective points of $f$ correspond to the inflection
points of
its characteristic curve.\par

{\bf Proof}. The projective points are the inflection points of ${\cal F}$,
i.e. the
second order contacts with its Legendrian tangent line. In projection $\pi$
the curve $c_f$ has the second order contact with a great circle, that is, an
inflection point.

Thus the Ghys theorem states that the characteristic curve of an orientation
preserving projective line diffeomorphism has at least 4 distinct inflection
points.

\proclaim Proposition 5.2. The characteristic curve $c_f$ is an embedded curve
transverse
to every meridian of the sphere (a great circle through the poles) and it
bisects the area of the sphere.\par

{\bf Proof}. First we show that $c_f$ is transverse to every meridian.
Meridians
correspond
to 2-dimensional subspaces in ${\bf R}^2_1 \times {\bf R}^2_2$ that
nontrivially intersect
both factors ${\bf R}^2_1$ and ${\bf R}^2_2$; such a space is automatically
Lagrangian since
it contains two linearly independent symplectically orthogonal vectors. The
curve $c_f$ lifts to ${\bf R}^2_1 \times {\bf R}^2_2$ as $\widetilde C
(\alpha)$
from the previous
section. Let $\widetilde C_1$ and $\widetilde C_2$ be its projections to
${\bf R}^2_1$ and ${\bf R}^2_2$.

If $c_f$ is tangent to a meridian then the plane generated by $\widetilde C$
and
$\dot {\widetilde C}$ nontrivially intersect ${\bf R}^2_1$ and ${\bf R}^2_2$
for some value of $\alpha$.
That is, a nontrivial linear combination of $\widetilde C_1$ and $\dot
{\widetilde C_1}$
(or $\widetilde C_2$ and $\dot {\widetilde C_2}$) vanishes. But this is
impossible because
both curves $\widetilde C_1$ and $\widetilde C_2$ are star shaped.

Next we show that $c_f$ makes exactly one turn around the sphere. Indeed $f$ is
isotopic to identity, therefore $c_f$ is homotopic to $c_{id}$ in the class of
immersed curves transverse to the meridians. The curve $c_{id}$ is the
equator, and the result follows.

\null

\vspace{6cm}

{}~\hspace{2,5cm}\special{illustration figghys_7.ps}

\vspace{0,2cm}

Inflection points of spherical curves is the subject of the {\it Tennis Ball
Theorem} (see \cite{arn2},\cite{arn5},\cite{arn6}): an embedded curve that
bisects the area of the
sphere has at least 4 distinct inflections. Therefore the tennis ball theorem
implies the theorem on 4 zeroes of the Schwarzian derivative.

{\bf Remark}. One may consider a Legendrian fibration $p: {\bf RP}^3
\rightarrow  S^2$
 which identifies
${\bf RP}^3$ with the space of oriented contact elements of the sphere. The
curve
$p({\cal F})$ is a front; generically it has cusps. The flattening points of
${\cal F}$
correspond to the vertices of the front $p({\cal F})$, i.e. its third order
contacts
with circles on the sphere. The derivative curve $c_f$ is the {\it derivative}
of the front $p({\cal F}$)- see \cite{arn6} for more details. Thus the Ghys
theorem
indeed concerns the vertices of a front on the sphere.

\vskip 0,5cm

{\bf Acknowledgement}. This paper is the result of a reflection over
the lecture by E. Ghys \cite{ghy}. We are grateful to him and
to V. Arnold, D. Fuchs and M. Kazarian for stimulating discussions.

\vfill\eject



\begin{thebibliography}{99}



\vskip 0,5cm

\bibitem{app}
{\sc P.~Appel},
``Sur les propri\'et\'es des cubiques gauches et le mouvement
h\'elicoidal d'un corps solide.'' Th\`ese, Paris 1876.

\bibitem{arn1}
{\sc V.I.~Arnold},
``Ramified covering $CP^2\to S^4$, hyperbolicity and projective
topology.'' Sib. Math. Journal {\bf 29} (1988) n.5, 36--47.

\bibitem{arn2}
{\sc V.I.~Arnold},
``Topological invariants of plane curves and caustics.'' AMS Univ. Lecture
Series N.5,
Providence, 1994.

\bibitem{arn3}
{\sc V.I.~Arnold},
``On the number of flattening points on space curves.'' Preprint Inst.
Mittag-Leffler, 1994.

\bibitem{arn4}
{\sc V.I.~Arnold},
``Remarks on the sextactic and other points of plane curves.''
Preprint, 1995.

\bibitem{arn5}
{\sc V.I.~Arnold},
``On topological properties of Legendre projections in contact geometry
of wave fronts.'' Algebra i Analysis (S. Petersbourg Math. J.) {\bf 6} (1994).

\bibitem{arn6}
{\sc V.I.~Arnold},
``Geometry of sperical curves and algebra of quaternions'', Uspekhi Mat. Nauk,
{\bf 50} (1995) n.1, 3-68 (in Russian).

\bibitem{ghy}
{\sc E.~Ghys}, ``Cercles osculateurs et g\'eom\'etrie lorentzienne''
Talk at the journ\'ee inaugurale du CMI, Marseille, February 1995.

\bibitem{gmo}
{\sc L.~Guieu, E.~Mourre, V.Yu.~Ovsienko},
``Theorem on six vertices of a plane curve via the Sturm theory.''Preprint CPT,
1995.

\bibitem{kne}
{\sc A.~Kneser},
H.Weber-Festschrift 170-180, Leipzig-Berlin, 1912.


\bibitem{muh}
{\sc S.~Mukhopadhyaya},
Bull.Calcutta Math.
Soc. 1, 1909.

\bibitem{hur}
{\sc A.~Hurwitz},
 Math. Ann. 57 (1903), 425-446.

\bibitem{stu}
{\sc J.C.F.~Sturm},
``M\'emoire sur les \'equations diff\'erentielles du second ordre.''
J. Math. Pures Appl. 1 (1836), 106-186.


\bibitem{tab1}
{\sc S.~Tabachnikov},
``Around four vertices.'' Russian Math. Surveys {\bf 45} (1990), n.1,
229--230.


\bibitem{tab2}
{\sc S.~Tabachnikov},
``The four vertex theorem revisited - Two variations on the old theme.''
 Preprint I. Newton Inst., 1995.



\end{thebibliography}
 \end{document}